\documentclass[12pt]{report}
\usepackage{amssymb,amsmath,makeidx,verbatim}
\newcommand{\C}{\mathbb{C}}

\newcommand{\be}{\begin{enumerate}}
	\newcommand{\ee}{\end{enumerate}}
\newcommand{\bq}{\begin{eqnarray*}}
	\newcommand{\eq}{\end{eqnarray*}}

\begin{document}
	\newcommand{\disp}{\displaystyle}
	\thispagestyle{empty}
	\begin{center}
		\textsc{A Paley-Wiener theorem for the Joint-Eigenspace Fourier transform on noncompact symmetric spaces.\\}
		\ \\
		\textsc{Olufemi O. Oyadare}\\
		\ \\
		Department of Mathematics,\\
		Obafemi Awolowo University,\\
		Ile-Ife, $220005,$ NIGERIA.\\
		\text{E-mail: \textit{femi\_oya@yahoo.com}}\\
	\end{center}
	\begin{quote}
		{\bf Abstract.} {\it This paper conducts a geometric analysis of the Joint-Eigenspace Fourier transform of the symmetric space of the noncompact type. Our study shows how the Poisson transform builds up the well-known Helgason Fourier transform for an analysis of the complete duality of the underlying symetric space. Among other results, we establish an inversion formula, a Plancherel formula and (as our main result) a Paley-Wiener theorem for the Joint-Eigenspace Fourier transform on any noncompact symmetric space.}
	\end{quote}
	
	\ \\
	\ \\
	\ \\
	\ \\
	\ \\
	\ \\
	\ \\
	\ \\
	\ \\
	\ \\ 
	\ \\
	\ \\
	\ \\
	\ \\
	$\overline{2020\; \textmd{Mathematics}}$ Subject Classification: $53C35, \;\; 43A90, \;\; 42A38$\\
	Keywords: Fourier Transform: Noncompact symmetric spaces: Harish-Chandra spherical transform: Helgason Fourier transform.\\
	
	\ \\
	{\bf \S1. Introduction.}
	
	\indent This paper considers another aspect of the duality of the symmetric spaces of the noncompact type in which the {\it Poisson transform} builds up the {\it Helgason Fourier transform} into the {\it Joint-Eigenspace Fourier transform,} thus leading to a complete duality of these spaces.
	
	\indent We consider a noncompact connected semisimple Lie groups $G$ with a maximal compact subgroup $K$ whose Iwasawa decomposition is $G=KAN.$ The structure theory of $G$ may be found in $[30.]$ and $[53.].$ The quotient space $X=G/K$ is then called the symmetric space of the noncompact type on which is defined the well-known Helgason Fourier transform. The Joint-Eigenspace Fourier transform is the Fourier transform on $X$ which factors into the Helgason Fourier transform on $X$ and the Poisson transform on $B=K/M$ (for the {\it centralizer} $M$ of $A$ in $K$). In this respect, we see the Joint-Eigenspace Fourier transform on $X$ as the completion of the Helgason Fourier transform on $X$ onto the joint-eigenspace of the underlying symmetric space. This crucial observation makes the fine structure of the joint-eigenspace representations to be readily available to the Joint-Eigenspace Fourier transform for further understanding of the symmetric spaces $X$ via the Joint-Eigenspace Fourier transform on $X.$
	
	\indent The Joint-Eigenspace Fourier transform was passively noticed and briefly mentioned in relation to the compact symmetric space in Helgason $[30.],$ p. $327,$ and referred to as \lq\lq {\it a genuine Fourier transform}" on $X$ whose range questions \lq\lq {\it have significant applications to differential equations."} This genuine Fourier transform was however never extracted nor studied out of the unfounded fear of having to deal with two genuine Fourier transforms on $X,$ when the Helgason Fourier transform on $X$ had already performed well enough to the satisfaction of many. The present paper shows the need not to entertain any such fear as the Joint-Eigenspace Fourier transform on $X$ is here shown to simply complete the journey of the Helgason Fourier transform on $X$ onto being a {\it full} joint-eigenspace Fourier transform on $X,$ by generalizing the Paley-Wiener theorem ($[23.]$) for the Helgason Fourier transform to the Paley-Wiener theorem for the Joint-Eigenspace Fourier transform.
	
	\indent The topic of Paley-Wiener type theorems via the Helgason Fourier transform on nocompact symmetric spaces and its applications into diverse fields of research is a vast and very active research area of mathematics since it was first considered by Paley and Wiener in the context of the Fourier transform on a space of functions of a complex variable. The more recent works of Alberti, G. S., Bartolucci, F., De Mari, F., De Vito, E. $[1.]$ developed a variation of the Helgason's theory of dual $G-$homogeneous pairs leading to the establishment of intertwining properties and inversion formulae for some Radon transforms and to an application to classical Radon and $X-$ray transforms in $\mathbb{R}^{3}.$ The paper of Arcozzi, N., Monguzzi, A., Peloso, M. M., Salvatori, M. $[2.]$ proved the Paley-Wiener type theorems for some spaces of holomorphic functions on the Siegel upper half space, which included the Hardy space, the weighted Bergmann spaces, the weighted Dirichlet spaces, the Drury-Arveson spaces and the Dirichlet spaces, as well as some structure theorems and applications.
	
	\indent Asta, D. M. $[6.]$ establishes the minimax rate of convergence for the kernel density estimate on symmetric spaces of noncompact type through the abstraction of the vector space operations on Euclidean spaces to symmetric spaces and a reparametrization of the Helgason Fourier transform on such spaces. The subject of reproducing kernels also gets the attention of geometers as Baranov, A., Belov, Y., Borichev, A. $[7.]$ generalizes a theorem of Young and Belov on the reproducing kernels in the Paley-Wiener space by showing the completeness of a biorthogonal system (to a complete system) and that of a minimal system of reproducing kernels. Indeed for other transforms,  Bardaro, C., Butzer, P. L., Mantellini, I., Schmeisser, G. $[8.]$ contains a new approach to the subject of Paley-Wiener theorem for the Mellin transforms based on the Riemannian surfaces of the logarithm in Mellin-Bernstein spaces. Knowing fully well that the theory of the Helgason Fourier transform essentially involves other notable transforms, Bartolucci, F., De Mari, F., Monti, M. $[11.]$ studied the unitarity of the Horocyclic Radon transform on symmetric spaces in the general case where the quasi regular representation of the group of isometries of the symmetric space is not irreducible nor is it square-integrable.
	
	\indent The reader may wish to consult recent works by Bernstein, S., Faustino, N. $[12.],$ Dalai, R. K., Ghosh, S., Srivastava, R. K. $[13.],$ Din, A. Y. $[14.],$ Dong, C.-P. $[15.],$ Eichinger, B., Woracek, H. $[23.],$ Gimperlein, H., Krotz, B., Kuit, J. J., Schlichtkrull, H. $[24.],$ Heins, M., Moucha, A., Roth, O. $[26.],$ Li, J., Lu, G., Yang, Q. $[31.],$ Lu, G., Yang, Q. $[32.],$ Pesenson, I. Z. $[39.],\;[40.],$ Picardello, M. A., Salvatori, M., Woess, W. $[41.],$ Sadiq, K., Tamasan, A. $[42.],$ Sherman, A. $[43.],$ Sonoda, S., Ishikawa, I., Ikeda, M. $[44.],$ Suzuki, K. $[45.]$ and Weiske, C., Yu, J., Zhang, G. $[56.]$ for various uses of the Helgason Fourier transform on different symmetric spaces and the necessity of the Paley-Wiener results in theory and applications.
	
	\indent The paper is organised as follows. \S2. contains some preliminaries on $G$ and its symmetric space $X=G/K.$ Here the necessity of having the Joint-Eigenspace Fourier transform on $X$ was motivated from the structure of the Harish-Chandra spherical transform on $G$ as an extension and a generalization of the Harish-Chandra spherical transform on $G,$ extended in order to encapsulate the Helgason Fourier transform of $X.$ We also prove an inversion and a Plancherel formulae for the Joint-Eigenspace Fourier transform and give its corresponding $c-$function in terms of the Harish-Chandra $c-$function. The image of the Helgason Fourier transform under the Poisson transform is the subject considered in some details in \S3. Here we justify the choice of the name \lq\lq {\it Joint-Eigenspace Fourier transform,}" verify that the Joint-Eigenspace Fourier transform is a Fourier transform on the whole of $X$, which factors into the Helgason Fourier transform on $X$ and the Poisson transform on $B=K/M,$ state a functional equation satisfied by this transform and prove the main result which gives a Paley-Wiener theorem for the Joint-Eigenspace Fourier transform on $X.$
	
	\ \\
	{\bf \S2. The Joint-Eigenspace Fourier transform on a symmetric space.}
	
	\indent Let $G$ be a connected semisimple Lie group with finite center, we denote its Lie algebra by $\mathfrak{g}$
	whose \textit{Cartan decomposition} is given as $\mathfrak{g} = \mathfrak{t}\oplus\mathfrak{p}.$ Denote by $\theta$ the \textit{Cartan involution} on $\mathfrak{g}$ whose collection of fixed points is $\mathfrak{t}.$
	We also denote by $K$ the analytic subgroup of $G$ with Lie
	algebra $\mathfrak{t}.$  $K$ is then a maximal compact subgroup of $G,$ $[50.].$
	Choose a maximal abelian subspace  $\mathfrak{a}$ of $\mathfrak{p}$ with algebraic
	dual $\mathfrak{a}^*$ and set $A =\exp \mathfrak{a}.$  For every $\lambda \in \mathfrak{a}^*$ put
	$$\mathfrak{g}_{\lambda} = \{X \in \mathfrak{g}: [H, X] =
	\lambda(H)X, \forall  H \in \mathfrak{a}\},$$ and call $\lambda$ a restricted
	root of $(\mathfrak{g},\mathfrak{a})$ whenever $\mathfrak{g}_{\lambda}\neq\{0\}$ ($[30.]$).
	Denote by $\mathfrak{a}'$ the open subset of $\mathfrak{a}$
	where all restricted roots are $\neq 0,$ and call its connected
	components the \textit{Weyl chambers.}  Let $\mathfrak{a}^+$ be one of the Weyl
	chambers, define the restricted root $\lambda$ positive whenever it
	is positive on $\mathfrak{a}^+$ and denote by $\triangle^+$ the set of all
	restricted positive roots. Members of $\triangle^+$ which form a basis for $\triangle$ and can not be written as a linear combination of other members of $\triangle^+$ are called \textit{simple,} $[50.].$ We then have the \textit{Iwasawa
		decomposition} $G = KAN$, where $N$ is the analytic subgroup of $G$
	corresponding to $\mathfrak{n} = \sum_{\lambda \in \triangle^+} \mathfrak{g}_{\lambda}$,
	and the \textit{polar decomposition} $G = K\cdot
	cl(A^+)\cdot K,$ with $A^+ = \exp \mathfrak{a}^+,$ and $cl(A^{+})$ denoting the closure of $A^{+}.$
	
	If we set $M = \{k \in K: Ad(k)H = H$, $H\in \mathfrak{a}\}$ and $M' = \{k
	\in K : Ad(k)\mathfrak{a} \subset \mathfrak{a}\}$ and call them the
	\textit{centralizer} and \textit{normalizer} of $\mathfrak{a}$ in $K,$ respectively, then (see $[27.]$);
	(i) $M$ and $M'$ are compact and have the same Lie algebra and
	(ii) the factor  $\mathfrak{w} = M'/M$ is a finite group called the \textit{Weyl
		group}. $\mathfrak{w}$ acts on $\mathfrak{a}^*_{\C}$ as a group of linear
	transformations by the requirement $$(s\lambda)(H) =
	\lambda(s^{-1}H),$$ $H \in \mathfrak{a}$, $s \in \mathfrak{w}$, $\lambda \in
	\mathfrak{a}^*_\mathbb{\C}$, the complexification of $\mathfrak{a}^*$.  We then have the
	\textit{Bruhat decomposition} $$G = \bigsqcup_{s\in \mathfrak{w}} B m_sB$$ where
	$B = MAN$ is a closed subgroup of $G$ and $m_s \in M'$ is the
	representative of $s$ (i.e., $s = m_sM$). The Weyl group invariant members of a space shall be denoted by the superscript $^{\mathfrak{w}}.$
	
	Some of the most important functions on $G$ are the \textit{spherical
		functions} which we now discuss as follows.  A non-zero continuous
	function $\varphi$ on $G$ shall be called a \textit{(zonal) spherical
		function} whenever $\varphi(e)=1,$ $\varphi \in C(G//K):=\{g\in
	C(G)$: $g(k_1 x k_2) = g(x)$, $k_1,k_2 \in K$, $x \in G\}$ and $f*\varphi
	= (f*\varphi)(e)\cdot \varphi$ for every $f \in C_c(G//K),$ where $(f \ast g)(x):=\int_{G}f(y)g(y^{-1}x)dy$ ($[9.]$). This
	leads to the existence of a homomorphism $\lambda :
	C_c(G//K)\rightarrow \C$ given as $\lambda(f) = (f*\varphi)(e)$.
	This definition is equivalent to the satisfaction of the functional relation $$\int_K\varphi(xky)dk = \varphi(x)\varphi(y),\;\;\;x,y\in G,\;[29.].$$
	
	\indent It has been shown by Harish-Chandra that spherical functions on $G$
	can be parametrized by members of $\mathfrak{a}^*_{\C}.$  Indeed every
	spherical function on $G$ is of the form $$\varphi_{\lambda}(x) = \int_Ke^{(i\lambda-p)H(xk)}dk,\; \lambda
	\in \mathfrak{a}^*_{\C},$$  $\rho =
	\frac{1}{2}\sum_{\lambda\in\triangle^+} m_{\lambda}\cdot\lambda,$ where
	$m_{\lambda}=dim (\mathfrak{g}_\lambda),$ and that $\varphi_{\lambda} =
	\varphi_{\mu}$ iff $\lambda = s\mu$ for some $s \in \mathfrak{w}$ ($[3.]$). Some of
	the well-known properties of spherical functions are $\varphi_{-\lambda}(x^{-1}) =
	\varphi_{\lambda}(x),$ $\varphi_{-\lambda}(x) =
	\bar{\varphi}_{\bar{\lambda}}(x),$ $\mid \varphi_{\lambda}(x) \mid\leq \varphi_{\Re\lambda}(x),$ $\mid \varphi_{\lambda}(x)\mid\leq \varphi_{i\Im\lambda}(x),$ $\varphi_{-i\rho}(x)=1,$ $\lambda \in \mathfrak{a}^*_{\C},$ while $\mid \varphi_{\lambda}(x) \mid\leq \varphi_{0}(x),\;\lambda \in i\mathfrak{a}^{*},\;x \in G,$ $[50.].$ Also if $\Omega$ is the \textit{Casimir operator} on $G$ then
	$$\Omega\varphi_{\lambda} = -(\langle\lambda,\lambda\rangle +
	\langle \rho, \rho\rangle)\varphi_{\lambda},$$ where $\lambda \in
	\mathfrak{a}^*_{\C}$ and $\langle\lambda,\mu\rangle
	:=tr(adH_{\lambda} \ adH_{\mu})$ for elements $H_{\lambda}$, $H_{\mu}
	\in {\mathfrak{a}}$ ($[55.]$). The elements $H_{\lambda}$, $H_{\mu}
	\in {\mathfrak{a}}$  are uniquely defined by the requirement that $\lambda
	(H)=tr(adH \ adH_{\lambda})$ and $\mu
	(H)=tr(adH \ adH_{\mu})$ for every $H \in {\mathfrak{a}}$ ([$27.$],
	Theorem $4.2$). Clearly $\Omega\varphi_0 = 0,$ $[35.].$
	
	\indent Due to a hint dropped by Dixmier in his discussion of some functional calculus,
	it is necessary to recall the notion of
	a \textit{`positive-definite'} function and then discuss the situation for
	positive-definite spherical functions.  We call a continuous function
	$f: G \rightarrow \C$ (algebraically) positive-definite whenever, for all
	$ x_1,\dots,x_m $ in $G$ and all $ \alpha_1,\dots,\alpha_m$ in $\C,$ we have $$\sum^m_{i,j=1}\alpha_i\bar{\alpha}_jf(x^{-1}_i x_j) \geq 0.$$  It
	can be shown $(cf.\;[27.])$ that $f(e) \geq 0$ and $|f(x)| \leq f(e)$ for every
	$x \in G$ implying that the space ${\wp}$ of all
	positive-definite spherical functions on $G$ is a subset of the
	space ${\mathfrak{F}}^{1}$ of all bounded spherical functions on $G,$ $[14.].$
	
	\indent We know, by the Helgason-Johnson theorem (See $[29.]$), that $${\mathfrak{F}}^{1}=
	\mathfrak{a}^*+iC_{\rho}$$ where $C_{\rho}$ is the convex hull of $\{s\rho: s \in
	\mathfrak{w}\}$ in $\mathfrak{a}^*$ ($[15.]$). Defining the \textit{involution} $f^*$ of $f$ as $f^*(x) =
	\overline{f(x^{-1})}$, it follows that $f = f^*$ for every $f \in
	{\wp}$, and if $\varphi_{\lambda} \in {\wp}$, then $\lambda$
	and $\bar{\lambda}$ are Weyl group conjugate, leading to a realization of $\wp$ as a subset of $\mathfrak{w} \setminus \mathfrak{a}^*_{\C}.$  ${\wp}$ becomes
	a locally compact Hausdorff space when endowed with the \textit{weak $^{*}-$topology} as a subset of $L^{\infty}(G)$.
	
	\indent We shall refer to $X=G/K$ as the symmetric space (of the groups $G$ of isometries) of the noncompact type. For $\lambda\in \mathfrak{a}^*_{\C},$ we denote $C_{c}(X)$ by $\mathcal{E}(X)$ and define the joint-eigenspace corresponding to $\lambda$ as $$\mathcal{E}_{\lambda}(X):=\{f\in \mathcal{E}(X): Df=\Gamma(D)(i\lambda)f\;\mbox{for }\;D\in \textbf{D}(X)\},$$ where $\textbf{D}(X)$ is the algebra of differential operators on $X$ that are invariant under all translations. See also $[51.]$ and $[54.].$ Given that $g\in G$ let $H(g),A(g)\in\mathfrak{a}$ be determined the Iwasawa decompositions of $g$ as $$g=k_{1} \exp H(g)n_{1}=n_{2} \exp A(g)k_{2}$$ in which $k_{1},k_{2}\in K$ and $n_{1},n_{2}\in N.$ It follows therefore that $A(g)=-H(g^{-1})$ as a function on $G$ and, as a function on $G/K\times K/M,$ we have that $$A(gK,kM)=A(k^{-1}g)=-H(g^{-1}k).$$
	
	\indent Let the notation be as specified above. The fundamental theorem of harmonic analysis on any such group $G$ is the explicit computation and/or characterization of a well-defined Fourier map on space $C_{c}^{\infty}(G)$ of test functions or on the Harish-Chandra Schwartz space $\mathcal{C}(G),$ (See $[3.],\;[4.],\;[5.]$), or on $\mathcal{C}^{p}(G)$ (for $0<p\leq 2$), (See $[46.],\;[48.],\;[49.]$). The two types being commonly considered are the scalar-valued Fourier transform map and the vector-valued Fourier transform map. The scalar-valued map is the well-known Harish-Chandra spherical Fourier transform on $G$ (also called the spherical transform) given as $f\mapsto \tilde{f}$ on $\mathcal{C}^{p}(G)$ and defined as $$\tilde{f}(\lambda)=\int_{G}f(g)\varphi_{-\lambda}(g)dg,$$ whenever the integral is absolutely convergent. Several authors have considered different aspect of the image of this map and the reader may see Oyadare $[35.]$ for an orientation. Notable results were established by Ehrenpreis and Mautner (${20.},$ $[21.],$ $[22.]$), Arthur ($[3.],\;[4.],\;[5.]$), Trombi-Varadarajan ($[25],$ $[50.]$), Trombi ($[46.],\;[48.],\;[49]$), Eguchi ($[16.],$ $[17.],$ $[18.],$ $[19.]$) and Oyadare ($[33.],\;[34.],\;[35.],\;[36.],\;[37.],\;[38.]$).
	
	\indent The major result concerning the scalar-valued Harish-Chandra Fourier transform $f\mapsto \tilde{f}$ on $G$ is the following which gives explicit for he map.
	
	\indent {\bf 2.1 Theorem} (Oyadare $[35.]$). The scalar-valued Harish-Chandra Fourier transform $f\mapsto \tilde{f}$ is a Fr$\acute{e}$chet linear topological Schwartz algebra isomorphism of $\mathcal{C}^{p}(G)$ with $$\{(\mathcal{H}\xi_{1})^{-1}\cdot h\cdot(\mathcal{H}\xi_{1})^{-1}:\;h\in\bar{\mathcal{Z}}({\mathfrak{F}}^{\epsilon})\}\subsetneq\mathcal{C}^{p}(\widehat{G}),$$ where $\bar{\mathcal{Z}}({\mathfrak{F}}^{\epsilon})$ is the Trombi-Varadarajan image $.\;\Box$
	
	\indent Notations are as in $[25.]$ and $[35.].$ Following the trend of results around the works of Trombi, it became also necessary to seek the fundamental theorem of harmonic analysis on $G$ via the operator-valued Harish-Chandra Fourier transform $f\mapsto \mathcal{F}(f)$ on $G.$ Employing a technique of L. Garding it bcame possible to characterize and to use a basis to explicitly compute its image. This was fully achieved by Oyadare $[38.]$
	
	\indent {\bf 2.2 Theorem} (Oyadare $[38.]$). The operator-valued Harish-Chandra Fourier transform $f\mapsto \mathcal{F}{f}$ is a Fr$\acute{e}$chet linear topological Schwartz algebra isomorphism of $\mathcal{C}^{p}(G)$ with the Schwartz Fr$\acute{e}$chet multiplication algebra $\mathcal{C}^{p}(\widehat{G})$ consisting of block matrices of the form $$((\mathfrak{F}_{B}(\breve{\alpha})_{(\gamma,m)}(\Lambda)\otimes\mathfrak{F}_{H}(\breve{\alpha})_{(\gamma,l)}(Q:\chi:\nu))_{\gamma\in F, (l,m)\in\mathbb{Z}^{2}})_{F\subset \hat{K},|F|<\infty}.\;\Box$$
	
	\indent Notations are as in $[10.]$ and $[38.].$ The case of $G=SL(2,\mathbb{R})$ is contained in $[10.]$ and $[52.],$ while the symmetric space was considered in $[19.].$
	
	\indent It is clear from the comparison of these major landmarks (Theorems $2.1$ and $2.2$) that the scalar-valued transform is an entry in the countable-matrix realization of the operator-valued transform, via the basis employed.
	
	\indent A closer consideration of the properties of the elementary spherical functions (used in the definition of the scalar-valued Harish-Chandra Fourier transform) shows that $$\varphi_{-\lambda}(g)=\varphi_{\lambda}(g^{-1}),$$ for all $g\in G,$ $\lambda\in \mathfrak{a}^*_{\C}.$ This equality allows us to re-write the scalar-valued Fourier transform formula as $$\tilde{f}(\lambda)=\int_{G}f(g)\varphi_{-\lambda}(g)dg$$ 
	$$\;\;\;\;\;\;\;\;\;\;=\int_{G}f(g)\varphi_{\lambda}(g^{-1})dg$$ $$\;\;\;\;\;\;\;\;\;\;\;\;=\int_{G}f(g)\varphi_{\lambda}(eg^{-1})dg$$ $$\;\;=(f*\varphi_{\lambda})(e),$$ where we have written $e$ for the identity element of $G$ and $*$ for the convolution of functions on $G.$ This computation shows that the scalar-valued Harish-Chandra Fourier transform is actually a convolution of every member of of the space $\mathcal{C}^{p}(G)$ with each of the elementary spherical functions on $G$ and that this convolution is then evaluated at $e.$ It would therefore be of general importance if we could consider a convolution of members of the space $\mathcal{C}^{p}(G)$ with each of the elementary spherical functions on $G$ (or with each of the {\it Eisenstein integral} on $G$) evaluated at any point of $G,$ which may not necessarily be at the identity element $e.$
	
	\indent Our first motivation is to therefore consider the map $$f\mapsto\mathcal{H}_{g}f,$$ for every $g\in G,$ $f\in \mathcal{C}^{p}(G),$ where $$(\mathcal{H}_{g}f)(\lambda)=(f*\varphi_{\lambda})(g).$$ It is clear that $(\mathcal{H}_{e}f)(\lambda)=(f*\varphi_{\lambda})(e)=\tilde{f}(\lambda).$ This is the Fourier transform on $G$ earlier considered (and named the {\it spherical-convolution transform}) in Oyadare $[33.]$ where some of its properties were discussed leading to a Placherel formula with the computation of its corresponding measure.
	
	\indent The natural map $\pi:G\rightarrow G/K$ also suggests an extension of the spherical-convolution Fourier map $f\mapsto\mathcal{H}_{g}f$ (for $g\in G$) to the map $f\mapsto\mathcal{H}_{x}f$ (for $x\in G/K$). This may be termed the group-to-symmetric-space extension of the spherical-convolution Fourier transform. It would also be of considerable importance to study this group-to-symmetric-space extension of the Fourier transform and relate it with the known Fourier transform on $G/K.$ This is our aim in this paper.
	
	\indent It is well known that a Fourier transform has already been defined and studied on $G/K,$ called the {\it Helgason Fourier transform.} Explicitly, if $f$ is a function defined on $X=G/K,$ then its Helgason Fourier transform $$f\mapsto \hat{f}$$ is given as $$\hat{f}(\lambda,b)=\int_{X}f(x)e^{(-i\lambda+\rho)(A(x,b))}dx$$ for all $\lambda\in \mathfrak{a}^*_{\C},$ $b\in B=K/M$ for which the integral exists $[29.].$ The Helgason Fourier transform is itself an extension (but not a generalization) of the Harish-Chandra Fourier transform as given in the following.
	
	\indent {\bf 2.3 Lemma} (Helgason $[29.],$ p. $224$). The Helgason Fourier transform is an extension of the Harish-Chandra (spherical) Fourier transform from $K-$invariant functions on $X$ to not-necessarily $K-$invariant functions on $X.$
	
	\indent {\bf Proof.} Consider any $K-$invariant function $f$ (on $X$) in the Helgason Fourier transform. We could then replace $f$ with its translations $f^{\tau(k)}$ and recall that $$A(\tau(k)x,b)=A(k\cdot x,b)=A(x,k^{-1}b),$$ $[27.].$ This makes the above formula of the Helgason Fourier transform to be independent of $b,$ which when integrated over $b$ and using the elementary fact that $\varphi_{\lambda}(gK)=\int_{K}e^{(i\lambda+\rho)(A(kg))}dk,$ leads to the formula of the Harish-Chandra Fourier transform$.\;\Box$
	
	\indent This means that the Helgason Fourier transform on $X$ reduces to the formula for the Harish-Chandra spherical Fourier transform when $K-$invariant functions are under consideration. This lemma then suggests explicitly that, for $K-$invariant functions $f$ on $X,$ the Helgason Fourier transform $f\mapsto \hat{f}$ is explicitly given as $$\hat{f}(\lambda,b)=\tilde{f}(\lambda)=(f*\varphi_{\lambda})(e)=(\mathcal{H}_{e}f)(\lambda).$$
	
	\indent Just as the spherical-convolution transform on $G$ generalizes the Harish-Chandra Fourier transform on $G$ (which has now been extended to the Helgason Fourier transform on $X$), one would like to seek a generalization of the Helgason Fourier transform on $X$ via the spherical-convolution transform's extension from $G$ to $X=G/K.$ This quest leads us to the consideration of the map $$f\mapsto \mathcal{H}_{x}f$$ for $x\in G/K=X$ defined as $$(\mathcal{H}_{x}f)(\lambda)=(f\times\varphi_{\lambda})(x),$$ where $\varphi_{\lambda}$ is the elementary spherical functions on $X$ and $\times$ is the convolution of functions on $X$ defined as $$(f_{1}\times f_{2})\circ\pi:=(f_{1}\circ\pi)*(f_{2}\circ\pi)$$ for the natural map $\pi:G\rightarrow G/K$ $[28.].$ The following is a cautionary observation in relation to the extensions mentioned above.
	
	\indent {\bf 2.4 Lemma.} The transformation map $f\mapsto \mathcal{H}_{x}f$ defined above for $x\in X,$ $f\in C(X),$ does not restrict to the Helgason Fourier transform on X.
	
	\indent {\bf Proof.} Denote the identity of $X$ by $\bar{e}.$ Then $$(\mathcal{H}_{\bar{e}}f)(\lambda)=(f\times\varphi_{\lambda})(\bar{e})=(f\times\varphi_{\lambda})(eK)=\int_{G}f(g\cdot eK)\varphi_{\lambda}(g^{-1}\cdot eK)dg,$$ which is not equal to $\hat{f}(\lambda,b).\;\Box$
	
	The above lemma does not erode the fact of Lemma $2.3$ and does not preclude the possibility of the transformation map $f\mapsto \mathcal{H}_{x}f$ coinciding with the Helgason Fourier transform on a special class of functions on $X,$ but it however reveals that we now have a completely different transform on $X$ at hand which we now christen as follows.
	
	\indent {\bf 2.5 Defnition.} Let $x\in X=G/K.$ We shall refer to the transformation map $f\mapsto \mathcal{H}_{x}f$ as the {\it Joint-Eigenspace Fourier transform} on $X$ which is given as $$(\mathcal{H}_{x}f)(\lambda)=(f\times\varphi_{\lambda})(x),$$ for $f\in C_{c}(X),$ $\lambda\in \mathfrak{a}^*_{\C}$ and the elementary spherical functions $\varphi_{\lambda}$ on $X.\;\Box$
	
	\indent The choice of this name will be made clear after Theorem $3.4$ ahead while more would still be said about the relationship between the Helgason Fourier transform on $X$ and the Joint-Eigenspace Fourier transform on $X$ after Theorem $2.8$ below. Specific examples and their detailed computations for the Joint-Eigenspace Fourier trasform on $X$ is to be found in $[8.],$ $[29.]$ and $[32.].$
	
	\indent For $0<p\leq 2$ and $x\in X,$ the image of the Joint-Eigenspace Fourier transform on $X$ is given as the set $$\mathcal{C}^{p}_{x}(\mathfrak{a}^*_{\C}):=\{\mathcal{H}_{x}f:\lambda\mapsto(\mathcal{H}_{x}f)(\lambda),f\in \mathcal{C}^{p}(X)\}\}.$$ It is thus very important to give explicit description of $\mathcal{C}^{p}_{x}(\mathfrak{a}^*_{\C}).$ The beauty inherent in the Joint-Eigenspace Fourier transform on $X$ (over and above the Helgason Fourier transform) is that the map $f\mapsto\mathcal{H}_{x}f$ could be considered both as a functions of $x$ in $X$ and as a function of $\lambda$ in $\mathfrak{a}^*_{\C}.$
	
	\indent The defining convolution of the Joint-Eigenspace Fourier transform on $X$ was briefly considered  as a function of $X$ by Helgason $[29.]$ (but not as a Fourier transform on $X$) leading to its integral representation in terms of the Helgason Fourier transform on $X$ given as follows.
	
	\indent {\bf 2.6 Lemma} (Helgason $[29],$ p. $225$). Given that $f\in \mathcal{D}(X):=C^{\infty}_{c}(X)$ and $\lambda\in\mathfrak{a}^*_{\C},$ then $$(\mathcal{H}_{x}f)(\lambda)=\int_{B}e^{(i\lambda+\rho)(A(x,b))}\hat{f}(\lambda,b)db.$$
	
	\indent {\bf Proof.} $(\mathcal{H}_{x}f)(\lambda)=(f\times\varphi_{\lambda})(x)=\int_{G}f(g\cdot\bar{e})\varphi_{\lambda}(g^{-1}\cdot x) dg$ as an integral over $G.$ Using the symmetry of the spherical function $\varphi_{\lambda}$ (as a function on $G$) given as $\varphi_{\lambda}(g^{-1}h)=\int_{K}e^{(i\lambda+\rho)(A(kh))}e^{(-i\lambda+\rho)(A(kg))}dk$ (Helgason $[29],$ p. $224$) and the identity $A(gK,kM)=A(k^{-1}g)$ (for $g\in G,$ $k\in K$) we have the result$.\;\Box$
	
	When further considered  as a function of $x,$ Helgason computed the Helgason Fourier transform of the Joint-Eigenspace Fourier transform showing the intimate interconnection between the two Fourier transforms on $X$ which is restated in the notation of Definition $2.5$ as follows.
	
	\indent {\bf 2.7 Lemma.} Let $f\in \mathcal{D}(X).$ Then $$\hat{(\mathcal{H}_{x}f)}(\lambda,b)=\hat{f}(\lambda,b)\cdot\hat{\varphi_{\mu}}(\lambda),$$ $\lambda,\mu\in\mathfrak{a}^*_{\C}.$
	
	\indent {\bf Proof.} We already know that $(\mathcal{H}_{x}f)(\mu)=(f\times\varphi_{\mu})(x),$ whose Helgason Fourier transform is given as Lemma $1.4$ in Helgason $[29],$ p. $226.\;\Box$
	
	It is clear that the elementary spherical function $\varphi_{\mu}$ are $K-$bi-invariant as functions on $X.$ In particular, each $\varphi_{\mu}$ is $K-$invariant, so that $\hat{\varphi_{\mu}}(\lambda)=\tilde{\varphi_{\mu}}(\lambda).$ The import of the result of Lemma $2.7$ is that the Helgason Fourier transform of the Joint-Eigenspace Fourier transform (of $f\in C_{c}(X),$ when the Joint-Eigenspace Fourier transform is considered as a function on $X=G/K$) is a non-zero constant multiple (by the constant $c=\tilde{\varphi_{\mu}}(\lambda)$) of the (oridinary) Helgason Fourier transform of $f.$ This connection would help to deduce results for the Joint-Eigenspace Fourier transform from results already established for the Helgason Fourier transform $$f\mapsto\hat{f}(\lambda,b)=\frac{1}{\tilde{\varphi_{\mu}}(\lambda)}\hat{(\mathcal{H}_{x}f)}(\lambda,b).$$
	
	\indent It appears that there was no much interest in any consideration of the Joint-Eigenspace Fourier transform on $X$ as a function of $\lambda\in\mathfrak{a}^*_{\C}$ in much similarity with the detailed treatments of both the Harish-Chandra spherical Fourier transform $f\mapsto\tilde{f}(\lambda)$ and the Helgason Fourier transform $f\mapsto\hat{f}(\lambda,b),$ as such a treatment of the Joint-Eigenspace Fourier transform on $X$ (as a map $f\mapsto(\mathcal{H}_{x}f)(\lambda)$ of $\lambda$) could have led to a better understanding of the joint-eigenspaces on $X.$ Our aim is to kick-start this study. It is however of interest to further note the striking similarity between the inversion formulae for both the Harish-Chandra spherical Fourier transform on $G$ (Barker $[10]$) and the Joint-Eigenspace Fourier transform on $X$ established as follows.
	
	\indent {\bf 2.8 Theorem.} For every $f\in \mathcal{D}(X)$ we always have $$f(x)=w^{-1}\int_{\mathfrak{a}^{*}}(\mathcal{H}_{x}f)(\lambda)|c(\lambda)|^{-2}d\lambda,$$ where $w$ is the order of the Weyl group and $c(\lambda)$ is the Harish-Chandra $c-$function.
	
	\indent {\bf Proof.} It has been shown in Theorem $1.3$ of Helgason $[29.],$ p. $225,$ that $$f(x)=w^{-1}\int_{\mathfrak{a}^{*}}\int_{B}e^{(i\lambda+\rho)(A(x,b))}\hat{f}(\lambda,b))|c(\lambda)|^{-2}d\lambda db.$$ Hence, $$f(x)=w^{-1}\int_{\mathfrak{a}^{*}}\int_{B}e^{(i\lambda+\rho)(A(x,b))}\hat{f}(\lambda,b))|c(\lambda)|^{-2}d\lambda db$$ $$=w^{-1}\int_{\mathfrak{a}^{*}}(\int_{B}e^{(i\lambda+\rho)(A(x,b))}\hat{f}(\lambda,b))db)|c(\lambda)|^{-2}d\lambda$$ $$=w^{-1}\int_{\mathfrak{a}^{*}}(\mathcal{H}_{x}f)(\lambda)|c(\lambda)|^{-2}d\lambda,$$ by the use of Lemma $2.6$ above$.\;\Box$
	
	\indent The inversion formula for the Joint-Eigenspace Fourier transform on $X,$ given in Theorem $2.8$ above, was curiously notes in $[29.],$ p. $327,$ as the noncompact analogue of the convergent expansion of an arbitrary function on a compact symmetric space and, with the perspective of the Harish-Chandra's inversion formula for the $K-$bi-invariant functions $$F(h)=\int_{K}f(ghkK)dk\;\;\;(f\in C_{c}(X)),$$ we have that $$F(h)=c\int_{\mathfrak{a}^{*}}\hat{F}(\lambda)\varphi_{\lambda}(h)|c(\lambda)|^{-2}d\lambda$$ into which we could insert $\hat{F}(\lambda)=(f\times\varphi_{\lambda})(gK)=(f\times\varphi_{\lambda})(x)=(\mathcal{H}_{x}f)(\lambda).$ Hence Theorem $2.8$ Our proof in Theorem $2.8$ above is however from the independent recognition of the Joint-Eigenspace Fourier transform as a Fourier transform on $X.$
	
	\indent This now makes it clear that, just as the Helgason Fourier transform is an extension of the Harish-Chandra spherical Fourier transform from $K-$invariant functions on $X$ to not necessarily $K-$invariant functions on $X$ (Recall Lemma $2.3$), the following is the two way bridge between the Helgason Fourier transform on $X$ and the Joint-Eigenspace Fourier transform on $X.$
	
	\indent {\bf 2.9 Theorem.} Let $f\in \mathcal{D}(X)$ and consider the $K-$bi-invariant functions $F(h)=\int_{K}f(ghkK)dk$ in $\textbf{D}(G).$ Then the Joint-Eigenspace Fourier transform on $X$ of any $f\in \mathcal{D}(X)$ is exactly the Helgason Fourier transform of $F$ (on $G$) as an extension of the Harish-Chandra spherical transform to $X.$ That is, $$(\mathcal{H}_{x}f)(\lambda)=\hat{F}(\lambda),\;\;\lambda\in\mathfrak{a}^{*}.\;\Box$$
	
	It is clear that the Helgason Fourier transform $\hat{F}(\lambda)$ in the above Theorem $2.9$ is essentially the extended Harish-Chandra spherical transform of the $K-$bi-invariant functions $F.$ This theorem gives a formal link between these two genuine Fourier transforms on a noncompact symmetric space and is the two-way correspondence between them.
	
	\indent Theorem $2.8$ directly leads to the Plancherel formula for the Joint-Eigenspace Fourier transform on $X$ for which we set $$\mathfrak{a}^{*}_{+}=\{\lambda\in\mathfrak{a}:A_{\lambda}\in\mathfrak{a}^{+}\},$$ where $A_{\lambda}$ is characterized by the Killing form requirement that $B(A_{\lambda},H)=\lambda(H)$ for all $H\in\mathfrak{a}.$
	
	\indent {\bf 2.10 Theorem.} The Joint-Eigenspace Fourier transform $f\mapsto(\mathcal{H}_{x}f)(\lambda)$ on $X$ extends to an isometry of $L^{2}(X)$ onto $L^{2}(\mathfrak{a}^{*}_{+}\times B),$ whose measure is given as $|\hat{\varphi_{\mu}}(\lambda)|^{-2}|c(\lambda)|^{-2}d\lambda db$ on $\mathfrak{a}^{*}_{+}\times B.$ Moreover, $$\int_{X}|f(x)|^{2}dx=\int_{\mathfrak{a}^{*}_{+}\times B}|\hat{(\mathcal{H}_{x}f)}(\lambda,b)|^{2}|\hat{\varphi_{\mu}}(\lambda)|^{-2}|c(\lambda)|^{-2}d\lambda db.$$
	
	\indent {\bf Proof.} The Helgason Fourier transform $f\mapsto\hat{f}(\lambda,b)$ extends to an isometry of $L^{2}(X)$ onto $L^{2}(\mathfrak{a}^{*}_{+}\times B),$ whose measure is given as $|c(\lambda)|^{-2}d\lambda db$ on $\mathfrak{a}^{*}_{+}\times B$ ($[29],$ p. $227$), from where it is known that $$\int_{X}|f(x)|^{2}dx=\int_{\mathfrak{a}^{*}_{+}\times B}|\hat{f}(\lambda,b)|^{2}|c(\lambda)|^{-2}d\lambda db.$$ Using the fact of Lemma $2.7$ above establishes the result$.\;\Box$
	
	\indent Theorem $2.10$ above gives the group-to-symmetric-space extension of the Plancherel formula for the spherical convolution transform on $G$ established in Oyadare $[33.]$ as Theorem $4.7.$ It is observed that we still have on $X$ (as we had on $G$ in $[33.]$) that $|\hat{\varphi_{\mu}}(\lambda)|^{-2}|c(\lambda)|^{-2}d\lambda db$ is the Plancherel measure of the Joint-Eigenspace Fourier transform on $\mathfrak{a}^{*}_{+}\times B.$ This observation confirms that the $c-$function that is associated with the inversion of the Joint-Eigenspace Fourier transform on $X$ is simply given as $\hat{\varphi_{\mu}}(\lambda)c(\lambda).$
	
	\ \\
	{\bf \S3. Poisson image of the Helgason Fourier transform.}
	
	\indent A horocycle in $X=G/K$ is defined to be any orbit $N^{'}\cdot x,$ where $x\in X$ and $N^{'}$ is any subgroup of $G=KAN$ conjugate to $N.$ We shall denote the collection of horocycles on $X$ as $\Xi$ and endow it with the differentiable structure of $G/MN.$ $\Xi$ would be seen as the dual space of $X$ ($[1.],$ $[11.]$) under a very general transform on $X$ defined as follows.
	
	\indent {\bf 3.1 Definition.} The {\it Radon transform} of a function $f$ on $X$ is defined as $\hat{f}(\xi)=\int_{\xi}f(x)ds(x),$ for all $\xi\in\Xi$ for which the integral exists.
	
	\indent The Radon transform is an injective map as explicitly stated below.
	
	\indent {\bf 3.2 Theorem} ($[29],$ p. $104$). If $f\in L^{1}(X),$ then $\hat{f}(\xi)$ exists for almost all $\xi\in\Xi$ and if $\hat{f}(\xi)=0$ for almost all $\xi\in\Xi$ then $f(x)=0$ for almost all $x\in X.\;\Box$
	
	\indent Now let $\xi(x,b)$ denote the horocycle passing through the point $x\in X$ with normal $b\in B=K/M.$ We shall denote by $A(x,b)\in\mathfrak{a}$ the composite metric from the origin to $\xi(x,b).$ For $\lambda\in\mathfrak{a}^{*}_{\mathbb{C}},$ $b\in B,$ the function $$e_{\lambda,b}:x\mapsto e_{\lambda,b}(x):=e^{(i\lambda+\rho)(A(x,b))}$$ is a joint-eigenfunction of $\textbf{D}(X)$ and belongs to the joint-eigenspace $\mathcal{E}(X).$ Considering $g\in G$ for which $x=gK\in G/K,$ the integral $$\varphi_{\lambda}(g)=\int_{B}e_{\lambda,b}(gK)=\int_{B}e^{(i\lambda+\rho)(A(gK,b))}$$ is an elementary spherical function on $G$ while $$u_{s}(x)=e_{-is\rho,b}(x)=e^{(s\rho+\rho)(A(x,b))}$$ is a harmonic function on $X$ for each $s\in W.$ We the harmonic function $u_{s=1}(x)=u_{1}(x)=e^{2\rho(A(x,b))}$ is known as the {\it Poisson kernel} ($[6.]$) and this informed the use of $e_{\lambda,b}(x)$ as a kernel to create a map $C(B)\rightarrow\mathcal{E}_{\lambda}(X)$ defined as follows.
	
	\indent {\bf 3.3 Definition.} The {\it Poisson transform} $P_{\lambda}$ of a function $F$ on $B$ is defined as $$P_{\lambda}(x)=\int_{B}e^{(i\lambda+\rho)(A(x,b))}F(b)db,$$ $x\in X.$
	
	\indent The Poisson transform $P_{\lambda}$ maps $C(B)$ into the joint-eigenspace $\mathcal{E}_{\lambda}(X).$ The first connection of the Poisson transform with the Joint-Eigenspace Fourier transform is via Lemma $2.6$ above, giving the Joint-Eigenspace Fourier transform on $X$ as a composition of the Poisson transform on the Helgason Fourier transform on $X$ given as follows.
	
	\indent {\bf 3.4 Theorem.} Given that $f\in \mathcal{D}(X)$ and $\lambda\in\mathfrak{a}^{*}_{\mathbb{C}}$ then $$(\mathcal{H}_{x}f)(\lambda)=(P_{\lambda}\hat{f}(\lambda,\cdot))(x),$$ where $\hat{f}(\lambda,\cdot)\in C(B),$ for each $\lambda.$ In particular, the Joint-Eigenspace Fourier transform on $X$ maps $C(X)$ into the joint-eihenspace $\mathcal{E}_{\lambda}(X).$
	
	\indent {\bf Proof.} It is already known, from Lemma $2.6,$ that $$(\mathcal{H}_{x}f)(\lambda)=\int_{B}e^{(i\lambda+\rho)(A(x,b))}\hat{f}(\lambda,b)db,$$ which by Definition $3.3,$ may be seen as the Poisson transform of the function $F=\hat{f}(\lambda,\cdot)$ in $C(B).$ That is, $$(\mathcal{H}_{x}f)(\lambda)=(P_{\lambda}\hat{f}(\lambda,\cdot))(x).\;\Box$$
	
	\indent In other words, the Joint-Eigenspace Fourier transform on $X$ is a Fourier transform on the whole of $X$ (as a lift of the Harish-Chandra spherical Fourier transform on $G$ to the symmetric space $X=G/K$) which does not restrict to the Helgason Fourier transform on $X,$ for all of the entire function space $\mathcal{D}(X).$ It may however be possible to extract a non-empty subspace of $\mathcal{D}(X)$ on which we have the identity $(\mathcal{H}_{x}f)(\lambda)\equiv\hat{f}(\lambda,b),$ for all $x\in X$ and for all $b\in B=K/M.$ Such a subspace of $\mathcal{D}(X),$ if it exists, would contain exactly those $f\in \mathcal{D}(X)$ on which $P_{\lambda}\equiv1,$ the identity map on $C(B).$ Let us therefore set $\mathcal{O}(X)$ to be defined as $$\mathcal{O}(X)=\{f\in\mathcal{D}(X):(\mathcal{H}_{x}f)(\lambda)\equiv\hat{f}(\lambda,b),\;x\in X,\;b\in B=K/M\},$$ whenever it exists ($[29.],$ $[51.]$).
	
	\indent A summary of the (inter-)relationships among the Helgason Fourier transform, the Poisson transform and the Joint-Eigenspace Fourier transform is presented by the following diagram: $$C(X)\rightarrow^{\hat{f}(\lambda,b)}\rightarrow C(B)\rightarrow^{P_{\lambda}}\rightarrow \mathcal{E}_{\lambda}(X),$$ whose composition is the Joint-Eigenspace Fourier transform on $X.$
	
	\indent The diagram justifies the choice of the name in Definition $2.5$ above. It explains how the Poisson transform completes the duality of the Helgason Fourier transform into a the joint-eigenspace $\mathcal{E}_{\lambda}(X)$ all of which are embodied in the Joint-Eigenspace Fourier transform on $X.$ This observation had also been made for the spherical convolution transform on $G$ in Oyadare $[33.].$ It is therefore clear from Theorem $3.4$ how the surjectivity of the Joint-Eigenspace Fourier transform on $X$ could be deduced from the properties of the Helgason Fourier transform and the Poisson transform.
	
	\indent If, in particular, $\lambda\in\mathfrak{a}^{*}_{\mathbb{C}}$ satisfies $Re(<i\lambda,\alpha>)>0$ for $\alpha\in\sum^{+},$ $b_{o}$ is the origin $eM$ in $B=K/M$ then, for $H\in\mathfrak{a}^{+}$ and $a_{t}=\exp tH,$ we have $$\lim_{t\rightarrow\infty}e^{(-i\lambda+\rho)(tH)}(\mathcal{H}_{a_{t}\cdot o}f)(\lambda)=c(\lambda)\hat{f}(\lambda,b_{o}).$$ Indeed, $$\lim_{t\rightarrow\infty}e^{(-i\lambda+\rho)(tH)}(\mathcal{H}_{a_{t}\cdot o}f)(\lambda)=\lim_{t\rightarrow\infty}e^{(-i\lambda+\rho)(tH)}\int_{B}e^{(i\lambda+\rho)(A(a_{t}\cdot o,b))}\hat{f}(\lambda,b)db$$ $$=c(\lambda)\hat{f}(\lambda,b_{o}),$$ by Theorem $3.16$ of $[29.],$ p. $120.$
	
	\indent {\bf 3.5 Lemma.} Given that $Re(<i\lambda,\alpha>)>0$ for $\alpha\in\sum^{+}$ and some $\lambda\in\mathfrak{a}^{*}_{\mathbb{C}},$ then, for $a_{t}=\exp tH,$ $H\in\mathfrak{a}^{+}$ and $b_{o}$ as the origin $eM$ in $B=K/M,$ we have that $$\lim_{t\rightarrow\infty}e^{(-i\lambda+\rho)(tH)}(\mathcal{H}_{a_{t}\cdot o}f)(\lambda)=c(\lambda)\hat{f}(\lambda,b_{o}).\;\Box$$
	
	\indent We equally have a functional equation for the Joint-Eigenspace Fourier transform on $X$ given as $$\int_{K}(\mathcal{H}_{gk\cdot x}f)(\lambda)dk=\varphi_{\lambda}(x)(\mathcal{H}_{g\cdot \bar{e}}f)(\lambda)$$ for $x\in X,$ $g\in G,$ $f\in C^{\infty}_{c}(X).$
	
	\indent In order to now establish a Paley-Wiener theorem for the Joint-Eigenspace Fourier transform on $X$ we shall examine closely both the Paley-Wiener theorem for the Helgason Fourier transform on $X$ and the bijectivity of the Poisson transform. To this end, we shall refer to a $C^{\infty}$ function $\psi(\lambda,b),$ $\lambda\in\mathfrak{a}^{*}_{\mathbb{C}},$ $b\in B,$ as being a {\it holomorphic function of uniform exponential type} and hence belonging to the space $\mathcal{H}^{R}(\mathfrak{a}^{*}\times B),$ if it is holomorphic in $\lambda$ and if there exists a constant $R\geq0$ such that for each $N\in\mathbb{Z}^{+},$ $$\sup_{\lambda\in\mathfrak{a}^{*}_{\mathbb{C}},b\in B}e^{-R|Im \lambda|}(1+|\lambda|)^{N}|\psi(\lambda,b)|<\infty,$$ where $Im \lambda$ is the imaginary part of $\lambda$ and $|\lambda|$ is its modulus. We set $$\mathcal{H}(\mathfrak{a}^{*}\times B)=\bigcup_{R>0}\mathcal{H}^{R}(\mathfrak{a}^{*}\times B)$$ and denote by $\mathcal{H}(\mathfrak{a}^{*}\times B)_{W}$ as members of $\mathcal{H}(\mathfrak{a}^{*}\times B)$ which are Weyl group invariant. Members of $\mathcal{H}(\mathfrak{a}^{*}\times B)_{W}$ may be viewed as being of the form $\psi(\lambda,b)=\psi_{\lambda}(b),$ so that $\mathcal{H}(\mathfrak{a}^{*}\times B)_{W}\subset C(B)=$ the domain of the Poisson transform, $P_{\lambda}.$ Thus the Weyl group invariance of members $\psi(\lambda,b)=\psi_{\lambda}(b)\in\mathcal{H}(\mathfrak{a}^{*}\times B)_{W}$ may be stated in terms of $P_{\lambda}$ as $P_{s\lambda}(\psi_{s\lambda})=P_{\lambda}(\psi_{\lambda})$ for each $s\in W.$ That is, $$\int_{B}e^{(is\lambda+\rho)(A(x,b))}\psi_{s\lambda}(b)db=\int_{B}e^{(i\lambda+\rho)(A(x,b))}\psi_{\lambda}(b)db$$ for each $s\in W.$ The following is the crucial Paley-Wiener theorem for the Helgason Fourier transform.
	
	\indent {\bf 3.6 Theorem} ($[29.],$ p. $270$). The Helgason Fourier transform $f\mapsto\hat{f}(\lambda,b)$ is a bijection of $C^{\infty}_{c}(X)$ onto $\mathcal{H}(\mathfrak{a}^{*}\times B)_{W}.$ Moreover, $\psi=\hat{f}$ is in $\mathcal{H}(\mathfrak{a}^{*}\times B)_{W}$ iff $supp(f)\subset Cl(B_{R}(0)).\;\Box$
	
	\indent The earlier diagram connecting the three transforms of this paper now becomes refined as $$C^{\infty}_{c}(X)\rightarrow^{\hat{f}(\lambda,b)}\rightarrow \mathcal{H}(\mathfrak{a}^{*}\times B)_{W}\rightarrow^{P_{\lambda}}\rightarrow \mathcal{E}_{\lambda}(X).$$
	
	\indent It is now very clear that it is not a coincidence that the Poisson transform is naturally included in the definition of the Paley-Wiener space for the Helgason Fourier transform on $X.$ Our contribution in this respect is to explore this fact to complement the journey of the Helgason Fourier transform on $X.$ We believe that this natural involvement of the Poisson transform in the image-construction of the Paley-Wiener space for the Helgason Fourier transform on $X$ is meant to complement the Helgason Fourier transform and to build it up into the status of a Joint-Eigenspace Fourier transform on $X.$ We now embark on the completion of the building up of Helgason Fourier transform on $X.$
	
	\indent Let $e^{\beta}$ be the $e_{s^{*}}-$function for $G_{\beta}/K_{\beta}$ in which $$e_{s^{*}}(\lambda)=const\cdot\prod_{\beta\in\sum^{+}_{o}}e^{\beta}(\lambda_{\beta}).$$ We know, from $[29],$ p. $269,$ that $P_{\lambda}$ is injective iff $e_{s^{*}}(\lambda)\neq0.$ On the surjectivity of $P_{\lambda}$ we note that if $\nu:\textbf{D}(G)\rightarrow\textbf{E}(X)$ denote the homomorphism induced by the action $f\mapsto f^{\tau(g)}$ of $G$ on $\mathcal{E}(X)$ and if $$\mathcal{E}^{\infty}_{\lambda}(X)=\{f\in\mathcal{E}_{\lambda}(X): (\nu(D)f)(x)=O(e^{A\;d(0,x)}),\forall D\in\textbf{D}(G)\;\mbox{and some}\;A>0\},$$ then for each $\lambda\in\mathfrak{a}^{*}_{\mathbb{C}}$ for which $e_{s^{*}}(\lambda)\neq0,$ we have that $$P_{\lambda}(\mathcal{E}(B))=\mathcal{E}^{\infty}_{\lambda}(X).$$
	
	\indent Since the $\mathcal{H}(\mathfrak{a}^{*}\times B)_{W}\subset C(B)$ we can now restrict $P_{\lambda}$ to $\mathcal{H}(\mathfrak{a}^{*}\times B)_{W}.$ Clearly $P_{\lambda}(\mathcal{H}(\mathfrak{a}^{*}\times B)_{W})\subset\mathcal{E}^{\infty}_{\lambda}(X)$ and if we define the Hilbert space $H^{\infty}_{\lambda}(X)$ as $$H^{\infty}_{\lambda}(X):=H_{\lambda}(X)\bigcap\mathcal{E}^{\infty}_{\lambda}(X),$$ from the Hilbert space $$H_{\lambda}(X):=\{h\in\mathcal{E}_{\lambda}(X):h(x)=\int_{B}e^{(i\lambda+\rho)(A(x,b))}F(b),\;F\in L^{2}(B)\}$$ ($[29.],$ p. $531$), then $$P_{\lambda}(\mathcal{H}(\mathfrak{a}^{*}\times B)_{W})=H^{\infty}_{\lambda}(X),$$ where $\lambda\in\mathfrak{a}^{*}_{\mathbb{C}}$ is simple. Thus a Paley-Wiener theorem for the Joint-Eigenspace Fourier transform is immediate and is as given in the following main result of the paper.
	
	\indent {\bf 3.7 Theorem.} Let $\lambda\in\mathfrak{a}^{*}_{\mathbb{C}}$ be simple. The Joint-Eigenspace Fourier transform on $X$ is a bijection of $C^{\infty}_{c}(X)$ onto $H^{\infty}_{\lambda}(X).$ Moreover, we have that $\psi(x)=(\mathcal{H}_{x}f)(\lambda)$ is in $H^{\infty}_{\lambda}(X)$ iff $supp(f)\subset Cl(B_{R}(0)).$
	
	\indent {\bf Proof.} The first part of this result is proven from the earlier computations when combined with Theorem $3.6.$ The second assertion is established as follows.
	
	$supp(f)\subset Cl(B_{R}(0))$ iff $\psi=\hat{f}$ is in $\mathcal{H}(\mathfrak{a}^{*}\times B)_{W}$ (Theorem $3.6$)\\
	$$\mbox{iff}\;\psi(\lambda,b)=\hat{f}(\lambda,b)\;\mbox{is in}\;\mathcal{H}(\mathfrak{a}^{*}\times B)_{W},\;b\in B,\;\mbox{for}\;\lambda\in\mathfrak{a}^{*}_{\mathbb{C}}\;\mbox{simple}$$
	$$\mbox{iff}\;\psi_{\lambda}(b)=\hat{f}(\lambda,b)\;\mbox{is in}\;\mathcal{H}(\mathfrak{a}^{*}\times B)_{W},\;b\in B,\;\mbox{for}\;\lambda\in\mathfrak{a}^{*}_{\mathbb{C}}\;\mbox{simple}\;$$
	$$\mbox{iff}\;\psi_{\lambda}(\cdot)=\hat{f}(\lambda,\cdot)\;\mbox{is in}\;\mathcal{H}(\mathfrak{a}^{*}\times B)_{W},\;\mbox{for}\;\lambda\in\mathfrak{a}^{*}_{\mathbb{C}}\;\mbox{simple}\;\;\;\;\;\;\;\;\;\;\;\;\;$$
	$$\mbox{iff}\;P_{\lambda}(\psi_{\lambda}(\cdot))=P_{\lambda}(\hat{f}(\lambda,\cdot))\;\mbox{is in}\;P_{\lambda}(\mathcal{H}(\mathfrak{a}^{*}\times B)_{W})=H^{\infty}_{\lambda}(X),\;\mbox{for}\;\lambda\in\mathfrak{a}^{*}_{\mathbb{C}}\;\mbox{simple}$$
	
	(which follows from the injectivity of $P_{\lambda}$ $[29.],$ p. $268-269$)
	$$\mbox{iff}\;P_{\lambda}(\psi_{\lambda}(\cdot))=(\mathcal{H}_{x}f)(\lambda)\;\mbox{is in}\;H^{\infty}_{\lambda}(X),\;\mbox{for}\;\lambda\in\mathfrak{a}^{*}_{\mathbb{C}}\;\mbox{simple}.\;\Box\;\;$$
	
	\indent It may be interesting to investigate Theorem $3.7$ for all $\lambda\in\mathfrak{a}^{*}_{\mathbb{C}}$ or on the Siegel upper-half plane. More properties of the Paley-Wiener space $H^{\infty}_{\lambda}(X)$ may be studied in the light of results in $[7.],$ $[12.],$ $[13.]$ and $[14.].$ We are now in a position to characterize the earlier defined $\mathcal{O}(X)-$subspace of the domain $\mathcal{D}(X)$ of the Joint-Eigenspace Fourier transform on $X.$
	
	\indent {\bf 3.8 Corollary.} Let $\lambda\in\mathfrak{a}^{*}_{\mathbb{C}}$ be simple. Then the subspace $\mathcal{O}(X)$ exists iff $\mathcal{H}(\mathfrak{a}^{*}\times B)_{W}=H^{\infty}_{\lambda}(X).$
	
	\indent {\bf Proof.} $\mathcal{O}(X)\neq\emptyset$ iff $(\mathcal{H}_{x}f)(\lambda)=\hat{f}(\lambda,b),$ for all $f\in\mathcal{O}(X),$ $x\in X,$ $b\in B$\\
	$$\;\;\;\;\;\;\;\;\mbox{iff}\;P_{\lambda}\equiv1,\mbox{the identity map on}\;\mathcal{H}(\mathfrak{a}^{*}\times B)_{W}$$
	$$\;\;\;\;\;\;\;\;\;\;\;\;\;\;\;\;\;\;\mbox{iff}\;H^{\infty}_{\lambda}(X)=P_{\lambda}(\mathcal{H}(\mathfrak{a}^{*}\times B)_{W})=I(\mathcal{H}(\mathfrak{a}^{*}\times B)_{W})=\mathcal{H}(\mathfrak{a}^{*}\times B)_{W}.\;\Box$$
	
	\indent The Paley-Wiener version of our diagram is now commutative and is finally given (for $\lambda\in\mathfrak{a}^{*}_{\mathbb{C}}$ simple) as 
 	 
    $$C^{\infty}_{c}(X)\rightarrow^{\hat{f}(\lambda,b)}\rightarrow \mathcal{H}(\mathfrak{a}^{*}\times B)_{W}\rightarrow^{P_{\lambda}}\rightarrow \mathcal{E}_{\lambda}(X).$$

    \ \\
	\ \\
	{\bf \S4. Conclusion.}
	
	The Poisson transform was shown here to be the appropriate complement to build up the Helgason Fourier on a noncompact symmetric space $X$ into the Joint-Eigenspace Fourier transform on $X.$ In specific terms, we have proved that the Joint-Eigenspace Fourier transform on $X$ is the Fourier transform on $X$ which factors into the Helgason Fourier transform (on $X=G/K$) and the Poisson transform (on $B=K/B$). The beautiful advantage of the Joint-Eigenspace Fourier transform on $X$ over and above the Helgason Fourier transform on $X$ is that the Joint-Eigenspace Fourier transform on $X$ naturally has its image in the joint-eigenspace $\mathcal{E}_{\lambda}(X)$ of $X$ whose fine differential structure and representations could now be fully developed (as against the development of only the joint-eigensubspaces of $\mathcal{E}_{\lambda}(X)$) and would further enhance the understanding of the Joint-Eigenspace Fourier transform on $X.$
	
	\indent The initial fear of having to deal with two genuine Fourier transforms on a noncompact symmetric space $X$ (each having distinct formulae, distinct inversion formulae, distinct Plancherel formulae, distinct $c-$functions and distinct Paley-Wiener theorems), is now allayed by factoring of the Joint-Eigenspace Fourier transform on $X$ as the Fourier transform on $X$ which factors into the Helgason Fourier transform (on $X=G/K$) and the Poisson transform (on $B=K/B$). This is to show that both the Joint-Eigenspace Fourier transform on $X$ and the Helgason Fourier transform on $X$ are genuine Fourier transforms on $X$ in which is seen that the Poisson transform on $B=K/M$ serves as the completion-transform to the Helgason Fourier transform in order to get it lifted to attain the beautiful status of the Joint-Eigenspace Fourier transform on $X.$ It is therefore appropriate to refer to the Helgason Fourier transform as the Weyl group-invariant Poisson-space Fourier transform on $X$ as against the present Joint-Eigenspace Fourier transform on $X.$
	
	\indent In explicit terms, the greatest advantage of the Joint-Eigenspace Fourier transform on $X$ above the Helgason Fourier transform on $X$ is that the Helgason Fourier transform on $X$ has to be (artificially) imported into the important analysis of differential equations of functions on $X$ in the consideration of the joint-eigenspace functions on $X$ while (on the other hand) the Joint-Eigenspace Fourier transform on $X$ naturally lives in the joint-eigenspace $\mathcal{E}_{\lambda}(X),$ $\lambda\in\mathfrak{a}^{*}$ simple, and naturally satisfies the differential equations $$D((\mathcal{H}_{x}f)(\lambda))=\Gamma(D)(i\lambda)(\mathcal{H}_{x}f)(\lambda),$$ for $D\in\textbf{D}(X),$ $x\in X$ ($[26.]$ and $[29.]$). It is now a matter of one's targeted destination of the needed Fourier transform on $X$ (whether the Weyl group-invariant Poisson-space $\mathcal{H}(\mathfrak{a}^{*}\times B)_{W}$ or the Joint-Eigensubspace $H^{\infty}_{\lambda}(X)$), which Fourier transform on $X$ (whether the Helgason Fourier transform or the Joint-Eigenspace Fourier transform, respectively) a researcher wishes to employ to achieve the aims of the research. It is however clear that the Joint-Eigenspace Fourier transform on $X$ does not just have significant contributions or applications to differential equations on $X$ but that the Joint-Eigenspace Fourier transform on $X$ naturally lives in the joint-eigenspace $\mathcal{E}_{\lambda}(X)$ and auromatically satisfies the aforementioned differential equations on $X.$ Our knowledge of the classical Euclidean case helps us to make bold to say that the Fourier transform on any noncompact symmetric space $X$ is to be that of the Joint-Eigenspace Fourier transform of this paper. This equally gives a complete version of the duality on such symmetric spaces.
	
	\indent An operator form of this theory may be conducted for the Joint-Eigenspace Fourier transform $f\mapsto(\mathcal{H}_{x}f)(\lambda),$ $x\in X,$ where $$(\mathcal{H}_{x}f)(\lambda):=(f\times \Phi_{\lambda,\delta})(x),$$ with $\delta\in K_{M},$ $\lambda\in\mathfrak{a}^{*}_{\mathbb{C}}$ and $\Phi_{\lambda,\delta}$ is the Eisenstein integral of $X$ ($[17.],$ $[47.]$). We shall consider the contribution of the Joint-Eigenspace Fourier transform to the theory of Eigenspace representation of functions and distributions on $X$ and the $L^{p}-$version of these results in another paper.
	\ \\
	\ \\
	\ \\
	{\bf   References.}
	\begin{description}
		
		\item [{[1.]}] Alberti, G. S., Bartolucci, F., De Mari, F., De Vito, E., Radon transform: Dual Pairs and Irreducible Representations, \textit{arXiv:$2002.01165v1.$} [math.FA], $2020.$
		
		\item [{[2.]}] Arcozzi, N., Monguzzi, A., Peloso, M. M., Salvatori, M., Paley-Wiener theorems on the Siegel upper half-space, \textit{arXiv:$1710.10079v2.$} [math.CV], $2024.$
			
		\item [{[3.]}] Arthur, J. G., \textit{Harmonic analysis of tempered distributions on semi-simple Lie groups of real rank one}, Ph.D. Dissertation, Yale University, $1970.$
		
		\item [{[4.]}] Arthur, J. G., \textit{Harmonic analysis of the Schwartz space of a reductive Lie group I,} mimeographed note, Yale University, Mathematics Department, New Haven, Conn.
		
		\item [{[5.]}] Arthur, J. G., \textit{Harmonic analysis of the Schwartz space of a reductive Lie group II,} mimeographed note, Yale University, Mathematics Department, New Haven, Conn.
		
		\item [{[6.]}] Asta, D. M., Lower bounds for kernel density estimation on symmetsric spaces, \textit{arXiv:$2403.10480v1.$} [math.ST], $2024.$
		
		\item [{[7.]}] Baranov, A., Belov, Y., Borichev, A., The Yong type theorem in weighted Fock spaces, \textit{arXiv:$1705.05778v1.$} [math.CV], $2017.$
		
		\item [{[8.]}] Bardaro, C., Butzer, P. L., Mantellini, I., Schmeisser, G., A fresh approach to the Paley-Wiener theorem for Mellin transforms and the Mellin-Hardy spaces, \textit{arXiv:$1706.00285v1.$} [math.FA], $2017.$
		
		\item [{[9.]}] Bargmann, V.,  Irreducible unitary representations of the Lorentz group, {\it Ann. of Math.} vol. {\bf 48}, $(1947),$ p. $568-640..$
		
		\item [{[10.]}] Barker, W. H.,  $L^{p}$ harmonic analysis of $SL(2,\mathbb{R}),$ \textit{American Mathematical Society Memoirs,} vol. \textbf{76} No. \textbf{393}, ($1988$).
		
		\item [{[11.]}] Bartolucci, F., De Mari, F., Monti, M., Unitarization of the Horocyclic Radon transform on symmetsric spaces, \textit{arXiv:$2108.04338v1.$} [math.RT], $2021.$
		
		\item [{[12.]}] Bernstein, S., Faustino, N., Paley-Wiener type theorems associated to Dirac operators of Riesz-Feller type, \textit{arXiv:$2405.04989v1.$} [math.CV], $2024.$
		
		\item [{[13.]}] Dalai, R. K., Ghosh, S., Srivastava, R. K., Quarternion Weyl transform and some uniqueness results, \textit{arXiv:$2110.00396v1.$} [math.FA], $2021.$
		
		\item [{[14.]}] Din, A. Y., A Paley-Wiener theorem for spherical $p-$adic spaces and Bernstein morphisms, \textit{arXiv:$2002.10063v2.$} [math.RT], $2020.$
		
		\item [{[15.]}] Dong, C.-P., On the Helgason-Johnson bound, \textit{arXiv:$2012.13474v3.$} [math.RT], $2021.$
		
		\item [{[16.]}] Eguchi, M., The Fourier Transform of the Schwartz space on a semisimple Lie group, \textit{Hiroshima Math. J.}, $vol.$ \textbf{4}, ($1974$), pp. $133-209.$
		
		\item [{[17.]}] Eguchi, M., Asymptotic expansions of Eisenstein integrals and Fourier transforms on symmetric spaces, \textit{J. Funct. Anal.} \textbf{34}, ($1979$), pp. $167 - 216.$
		
		\item [{[18.]}] Eguchi, M., Some properties of Fourier transform on Riemannian symmetric spaces, \textit{Lecture on Harmonic Analysis on Lie Groups and related Topics,} (T. Hirai and G. Schiffmann (eds.)) Lectures in Mathematics, Kyoto University, No. \textbf{4}) pp. $9 - 43.$
		\item [{[19.]}] Eguchi, M. and Kowata, A., On the Fourier transform of rapidly decreasing function of $L^{p}$ type on a symmetric space, \textit{Hiroshima Math. J.} $vol.$ \textbf{6}, ($1976$), pp. $143 - 158.$
		
		\item [{[20.]}] Ehrenpreis, L. and Mautner, F. I., Some properties of the Fourier transform on semisimple Lie groups, I, \textit{Ann. Math.}, $vol.$ $\textbf{61}$ ($1955$), pp. 406-439;
		
		\item [{[21.]}] Ehrenpreis, L. and Mautner, F. I., Some properties of the Fourier transform on semisimple Lie groups, II, \textit{Trans. Amer. Math. Soc.}, $vol.$ $\textbf{84}$ ($1957$), pp. $1-55;$
		
		\item [{[22.]}] Ehrenpreis, L. and Mautner, F. I., Some properties of the Fourier transform on semisimple Lie groups, III, \textit{Trans. Amer. Math. Soc.}, $vol.$ $\textbf{90}$ ($1959$), pp. $431-484.$
		
		\item [{[23.]}] Eichinger, B., Woracek, H., Homogeneous spaces of entire functions, \textit{arXiv:$2407.04979v1.$} [math.CV], $2024.$
		
		\item [{[24.]}] Gimperlein, H., Krotz, B., Kuit, J. J., Schlichtkrull, H., A Paley-Wiener theorem for Harish-Chandra Modules,\\ \textit{arXiv:$2010.04464v3.$} [math.RT], $2021.$
		
		\item [{[25.]}] Gangolli, R. and Varadarajan, V. S., \textit{Harmonic analysis of spherical functions on real reductive groups,} Ergebnisse der Mathematik und iher Genzgebiete, $vol.$ {\bf 101}, Springer-Verlag, Berlin-Heidelberg. $1988.$
		
		\item [{[26.]}] Heins, M., Moucha, A., Roth, O., Spectral theory of the invariant Laplacian on the disk and the sphere: a complex analysis approach, \textit{arXiv:$2312.12900v2.$} [math.CV], $2023.$
		
		\item [{[27.]}] Helgason, S., \textit{Differential geometry and symmetric spaces,} Academic Press, New York, $1962.$
		
		\item [{[28.]}] Helgason, S., A duality for symmetric spaces with applications to group representations, \textit{Advances in Mathematics,} $vol.$ $\textbf{5}$ ($1970$), pp. $1-154.$
		
		\item [{[29.]}] Helgason, S., \textit{Geometric analysis on symmetric spaces,} Mathematical Surveys and Monographs, Providence, Rhode Island, $1994.$
		
		\item [{[30.]}] Knapp, A.W., \textit{Representation theory of semisimple groups; An overview based on examples,} Princeton University Press, Princeton, New Jersey. $1986.$
		
		\item [{[31.]}] Li, J., Lu, G., Yang, Q., Higher order Brezis-Nirenberg problem on hyperbolic spaces: Existence, Non-existence and Symmetry of solutions, \textit{arXiv:$2107.03967v1.$} [math.AP], $2021.$
		
		\item [{[32.]}] Lu, G., Yang, Q., Sharp Hardy-Sobolev-Maz'ya, Adams and Hardy-Adams inequalities on the Siegel domains and complex hyperbolic spaces, \textit{arXiv:$2106.02103v1.$} [math.CA], $2021.$
		
		\item [{[33.]}] Oyadare, O. O., On harmonic analysis of spherical convolutions on semisimple Lie groups, \textit{Theoretical Mathematics and Applications}, $vol.$ $\textbf{5},$ no.: {\bf 3}. ($2015$), pp. $19-36.$
		
		\item [{[34.]}] Oyadare, O. O., Series analysis and Schwartz algebras of spherical convolutions on semisimple Lie groups \textit{Algebras, Groups and Geometries} \textbf{40}(1), $(2024),$ p. $41-59.$ See also arXiv.$1706.09045$ [math.RT].
		
		\item [{[35.]}] Oyadare, O. O.,  Non-spherical Harish-Chandra Fourier transforms on real reductive groups, \textit{J. Fourier Anal. Appl.} \textbf{28}, $15$ $(2022).$\\ http://doi.org/10.1007/s00041-09906-w
		
		\item [{[36.]}] Oyadare, O. O.,  The full Bochner theorem on real reductive groups, \textit{Algebras, Groups and Geometries} \textbf{39}, $(2023),$ p. $207-220.$
		
		\item [{[37.]}] Oyadare, O. O.,  Functional analysis of canonical wave-packets on real reductive groups, \textit{arXiv:$1912.07542v1.$} [math.FA], $13$ Dec. $2019.$
		
		\item [{[38.]}] Oyadare, O. O.,  On the operator-valued Fourier transform of the Harish-Chandra Schwartz Algebra, \textit{arXiv:$2407.20755v1.$} [math.RT], $2024.$
		
		\item [{[39.]}] Pesenson, I. Z., Sobolev, Besov and Paley-Wiener vectors in Banach and Hilbert spaces, \textit{arXiv:$1708.07416v2.$} [math.FA], $2023.$
		
		\item [{[40.]}] Pesenson, I. Z., Besov and Paley-Wiener spaces, moduli of continuity and Hardy-Steklov opeartors associated with th group \lq\lq ax+b ", \textit{arXiv:$2401.16734v1.$} [math.FA], $2024.$
		
		\item [{[41.]}] Picardello, M. A., Salvatori, M., Woess, W., Polyharmonic potential theory on the Poincare disk, \textit{arXiv:$2312.05806v1.$} [math.FA], $2023.$
		
		\item [{[42.]}] Sadiq, K., Tamasan, A., On the range of the Planar $X-$ray transform on the Fourier lattice of the Torus, \textit{arXiv:$2201.10926v2.$} [math.AP], $2022.$
		
		\item [{[43.]}] Sherman, A., The Cartan-Helgason theorem for supersymmestric spaces: spherical weughts for Kac-Moody superalgebras, \textit{arXiv:$2403.19145v1.$} [math.RT], $2024.$
		
		\item [{[44.]}] Sonoda, S., Ishikawa, I., Ikeda, M., Fully-connected network on noncompact symetric spaces and ridgelet transform based on Helgason-Fourier analysis, \textit{arXiv:$2203.01631v2.$} [cs.LG], $2022.$
		
		\item [{[45.]}] Suzuki, K., Trace Paley-Wiener theorem for Braverman-Kazhdan's asymptotic Hecke Algebra, \textit{arXiv:$2407.02752v1.$} [math.RT], $2024.$
		
		\item [{[46.]}] Trombi, P. C., Spherical transforms on symmetric spaces of rank one (or Fourier analysis on semisimple Lie groups of split rank one), \textit{Thesis, University of Illinios} ($1970$).
		
		\item [{[47.]}] Trombi, P. C., On Harish-Chandra's theory of the Eisenstein integral for real semisimple Lie groups. \textit{University of Chicago Lecture Notes in Representation Theory,} ($1978$), pp. $287$-$350.$
		
		\item [{[48.]}] Trombi, P. C., Harmonic analysis of $\mathcal{C}^{p}(G:F)\;(1\leq p<2)$ \textit{J. Funct. Anal.,} $vol.$ {\bf 40}. ($1981$), pp. $84$-$125.$
		
		\item [{[49.]}] Trombi, P. C., Invariant harmonic analysis on split rank one groups with applications. \textit{Pacific J. Math.} $vol.$ \textbf{101}. no.: \textbf{1}.($1982$), pp. $223 - 245.$
		
		\item [{[50.]}] Trombi, P. C. and Varadarajan, V. S., Spherical transforms on semisimple Lie groups, \textit{Ann. Math.,} $vol.$ {\bf 94}. ($1971$), pp. $246$-$303.$
		
		\item [{[51.]}] Varadarajan, V. S., Eigenfunction expansions on semisimple Lie groups, in \textit{Harmonic Analysis and Group Representation}, (A. Fig$\grave{a}$  Talamanca (ed.)) (Lectures given at the $1980$ Summer School of the \textit{Centro Internazionale Matematico Estivo (CIME)} Cortona (Arezzo), Italy, June $24$ - July $9.$ vol. \textbf{82}) Springer-Verlag, Berlin-Heidelberg. $2010,$ pp. $351-422.$
		
		\item [{[52.]}] Varadarajan, V. S., \textit{An introduction to harmonic analysis on semisimple Lie groups,} Cambridge Studies in Advanced Mathematics, \textbf{161},  Cambridge University Press, $1989.$
		
		\item [{[53.]}] Varadarajan, V. S., Harmonic analysis on real reductive reductive groups, \textit{Lecture Notes in Mathematics,} \textbf{576}, Springer Verlag, $1977.$
		
		\item [{[54.]}] Wallach, N.,  {\it Harmonic analysis in homogeneous spaces,} Dekker, New York, $1973.$
		
		\item [{[55.]}] Warner, G., {\it Harmonic analysis on semisimple Lie groups, I.} Springer-Verlag, New York, $1972.$
		
		\item [{[56.]}] Weiske, C., Yu, J., Zhang, G., Cartan-Helgason theorem for Quarternionic symmetric and twistor spaces, \textit{arXiv:$2306.15090v2.$} [math.FA], $2023.$
	\end{description}
\end{document}